\documentclass[11pt]{amsart}

\usepackage{amsmath}
\usepackage{amscd}
\usepackage{amssymb}
\usepackage{latexsym}
\usepackage{mathabx}

\usepackage[arrow,curve,matrix,tips,2cell]{xy}
  \SelectTips{eu}{10} \UseTips
  \UseAllTwocells

\usepackage{calrsfs}
\usepackage[T1]{fontenc} 
\usepackage{textcomp}
\usepackage{times}
 \usepackage[scaled=0.92]{helvet}

\usepackage{geometry}
\newtheorem{theorem}{Theorem}[section]
\newtheorem{lemma}[theorem]{Lemma}
\newtheorem{corollary}[theorem]{Corollary}
\newtheorem{proposition}[theorem]{Proposition}
\theoremstyle{definition}
\newtheorem{remark}[theorem]{Remark}
\newtheorem{definition}[theorem]{Definition}
\newtheorem{notation}[theorem]{Notation}

\newcommand{\bgl}{\begin{equation}}
\newcommand{\egl}{\end{equation}}
\newcommand{\bgln}{\begin{eqnarray}} 
\newcommand{\egln}{\end{eqnarray}}
\newcommand{\bglnoz}{\begin{eqnarray*}} 
\newcommand{\eglnoz}{\end{eqnarray*}}
\newcommand{\btheo}{\begin{theorem}}
\newcommand{\etheo}{\end{theorem}}
\newcommand{\blemma}{\begin{lemma}}
\newcommand{\elemma}{\end{lemma}}
\newcommand{\bproof}{\begin{proof}}
\newcommand{\eproof}{\end{proof}}
\newcommand{\bdefin}{\begin{definition}}
\newcommand{\edefin}{\end{definition}}
\newcommand{\bprop}{\begin{proposition}}
\newcommand{\eprop}{\end{proposition}}
\newcommand{\bcor}{\begin{corollary}}
\newcommand{\ecor}{\end{corollary}}
\newcommand{\bgloz}{\begin{equation*}} 
\newcommand{\egloz}{\end{equation*}}
\newcommand{\bremark}{\begin{remark}}
\newcommand{\eremark}{\end{remark}}

\newcommand{\cO}{\mathcal O}
\newcommand{\cP}{\mathcal P}

\def\Az{\mathbb{A}}
\def\Cz{\mathbb{C}}
\def\Kz{\mathbb{K}}
\def\Lz{\mathbb{L}}

\def\Nz{\mathbb{N}}
\def\Qz{\mathbb{Q}}
\def\Rz{\mathbb{R}}

\def\Q{\mathbb{Q}}
\def\Zz{\mathbb{Z}}

\newcommand{\mfm}{\mathfrak m}
\newcommand{\mfn}{\mathfrak n}
\newcommand{\mfp}{\mathfrak p}
\newcommand{\mfq}{\mathfrak q}
\newcommand{\mfa}{\mathfrak a}
\newcommand{\mfb}{\mathfrak b}
\newcommand{\mfc}{\mathfrak c}
\newcommand{\mfr}{\mathfrak r}

\newcommand{\ma}{\mapsto}

\newcommand{\Aut}{{\rm Aut}\,}
\newcommand{\id}{{\rm id}}
 
\newcommand{\defeq}{\mathrel{:=}} 
\newcommand{\dop}{\text{: }} 

\newcommand{\rec}{{\rm rec}}
\newcommand{\lge}{\left\{} 
\newcommand{\rge}{\right\}}

\newcommand{\gekl}[1]{\lge #1 \rge}
 
\newcommand{\menge}[2]{\gekl{ #1 \dop #2 }}
\newcommand{\ab}{\mbox{{\tiny \textup{ab}}}}
\newcommand{\Gab}{G^{\ab}}

\newcommand{\hpsi}{\widecheck{\psi}}
\newcommand{\Hab}[1]{\widecheck{G}^{\mathrm{ab}}_{#1}}
\newcommand{\Gal}{\mathrm{Gal}}
\newcommand{\an}[1]{``#1''} 
\newcommand{\ti}{\tilde}
\newcommand{\lori}{\longrightarrow}
\newcommand{\Rarr}{\Rightarrow} 
\newcommand{\Larr}{\Leftarrow} 
\newcommand{\LRarr}{\Leftrightarrow} 
\newcommand{\lspan}{{\rm span}}
\newcommand{\plim}{\varprojlim}
\usepackage{rotating}
\newcommand{\acts}{\;\rotatebox[origin=c]{-90}{$\circlearrowright$}\;}

\begin{document}

\date{\today\ (version 1.0)} 
\title[Dynamics in the abelianized Galois group]{Reconstructing global fields from\\[1mm] dynamics in the abelianized Galois group}
\thanks{This paper supersedes \cite{CM}, of which it is the second (and final) part, dealing with the dynamical systems aspects of the theory. The first part \cite{CM1} dealt with physics aspects related to partition functions, and \cite{CLMSS} is a companion paper, only containing number theoretical results. 
Part of this work was done whilst the first two authors enjoyed the hospitality of the University of Warwick (special thanks to
Richard Sharp for making it possible). } 
\author[G.~Cornelissen]{Gunther Cornelissen}
\address{Mathematisch Instituut, Universiteit Utrecht, Postbus 80.010, 3508 TA Utrecht, Nederland}
\email{g.cornelissen@uu.nl, h.j.smit@uu.nl}
\author[X.~Li]{Xin Li}
\address{School of Mathematical Sciences, Queen Mary University of London, Mile End Road, London E1 4NS, United Kingdom}
\email{xin.li@qmul.ac.uk}
\author[M.~Marcolli]{Matilde Marcolli}
\address{Mathematics Department, Mail Code 253-37, Caltech, 1200 E.\ California Blvd.\ Pasadena, CA 91125, USA}
\email{matilde@caltech.edu}
\author[H.~Smit]{Harry Smit}


\begin{abstract} \noindent We study a dynamical system induced by the Artin reciprocity map for a global field. We translate the conjugacy of such dynamical systems into various arithmetical properties that are equivalent to field isomorphism, relating it to anabelian geometry. 
\end{abstract}

\maketitle

\section{Introduction}

In this paper, we look at the main constructions of class field theory from the point of view of topological dynamical systems. If $\Kz$ is a global field, the finite ideles  $\Az_{\Kz,f}^*$ act on the abelianized absolute Galois group $\Gab_{\Kz}$ through multiplication with their image under the reciprocity map. Hence the group of fractional ideals (finite ideles modulo idelic units $\widehat{\cO}_{\Kz}^*$), and also the monoid $I_{\Kz}$ of integral ideals, act on a topological monoid $X_{\Kz}$, defined as the quotient of $\Gab_{\Kz} \times \widehat{\cO}_{\Kz}$ by the subgroup $\{(\rec(u)^{-1}, u)| u \in \widehat{\cO}_{\Kz}^*\}$, where $\widehat{\cO}_{\Kz}$ is the set of integral adeles (cf. Definition~\ref{def:X_K}). It is the action $I_{\Kz} \acts X_{\Kz}$ that we study. Our main result says that natural equivalences of such dynamical systems for different fields (orbit equivalence, topological conjugacy) have a purely number theoretical meaning, called a ``reciprocity isomorphism". Referring to the companion paper \cite{CLMSS}, this is equivalent to equality of Dirichlet $L$-series, and isomorphism of fields. The main result is Theorem \ref{THM} below. 

\medskip

We believe our study sheds some new light on the \emph{anabelian} geometry of global fields, as we now briefly explain. The anabelian philosophy is, in the words of Gro\-then\-dieck (\cite{Gro}, footnote (3)) ``a construction which pretends to ignore [\dots] the algebraic equations which traditionally serve to describe schemes, [\dots] to be able to hope to reconstitute a scheme [\dots] from [\dots] a purely topological invariant [\dots]''.  
In the zero-dimensional case, the fundamental group plays no r\^ole, only the absolute Galois group, and we arrive at the famous result of Neukirch and Uchida  that the existence of an isomorphism of absolute Galois groups of global fields implies the existence of an isomorphism of those fields (In \cite{NeukirchInv}, Neukirch proved this for number fields that are Galois over $\Q$; in \cite{U3}, Uchida proved it for general number fields and in \cite{U}, he proved it for function fields,  cf.\ also \cite{Neu3} 12.2, Ikeda \cite{Ikeda} and unpublished work of Iwasawa).

Individually, neither the Dedekind zeta function (i.e., the Dirichlet $L$-series for the trivial character), nor the adele ring, nor the abelianized Galois groups of a global field determine that field uniquely up to isomorphism:
\begin{itemize}
\item $\Q(\sqrt[8]{3})$ and $\Q(\sqrt[8]{3 \cdot 2^4})$ have the same zeta function, but non-isomorphic adele rings (\cite{G}, \cite{Komatsu}, \cite{Perlis1});
\item $\Q(\sqrt[8]{2 \cdot 9})$ and $\Q(\sqrt[8]{2^5 \cdot 9})$ have isomorphic adele rings (\cite{Komatsu2});
\item $\Q(\sqrt{-3})$ and $\Q(\sqrt{-7})$ have isomorphic abelianized Galois groups ( \cite{AS}, \cite{Onabe}, \cite{Kubota});
\item Equality of zeta functions for global function fields is the same as isomorphism of their adele rings, which is the same as the corresponding curves having isogenous Jacobians (\cite{Tate}, \cite{Turner}). 
\end{itemize}

The dynamical system $I_{\Kz} \acts X_{\Kz}$ that we consider (which, after all, \emph{is} a topological space with a monoid action) can be considered as some kind of substitute for the absolute Galois group. It involves the abelianized Galois group and the adeles, but not the absolute Galois group. In this sense, it is ``not anabelian''.  Our main result in this paper and \cite{CLMSS} says that consideration of this system cancels out the ``defects'' of the above individual examples in exactly the right way: it determines the field up to isomorphism. In the spirit of Kronecker's programme, one wants to characterize a number field by structure that is ``internal'' to it (i.e., not using extensions of the field): this is the case for the dynamical system, since class field theory realizes Kronecker's programme for abelian extensions. On the other hand, anabelian geometry characterizes a number field by its absolute Galois group, an object whose ``internal'' understanding remains largely elusive and belongs to the Langlands programme.

\medskip

We end this introduction with a brief discussion of related work and open problems. For the field of rational numbers, the dynamical system $I_{\Qz} \acts X_{\Qz}$ was introduced and studied by Bost and Connes in \cite{BC} (in which $\Nz \acts  \hat \Zz$ is the action given by multiplication by positive integers in the profinite integers). The construction for a general field was given by Ha and Paugam in \cite{HP} and clarified in \cite{LLN} and \cite{RN}. All of those references take the point of view of noncommutative geometry (specifically, $C^*$-algebras and quantum statistical mechanical systems). We describe the link between the noncommutative point of view and our work by considering algebraic crossed product algebras in Section \ref{NC}. In \cite{Yal}, Prop.\ 8.2, the space $X_{\Kz}$ was identified with the so-called Deligne-Ribet monoid from \cite{D-R}. 

In \cite{CM1}, two of us have shown that such quantum statistical mechanical systems are isomorphic if and only if the corresponding number fields are isomorphic, \emph{if} the number fields are Galois over $\Qz$ (analogue of the original theorem of Neukirch), essentially because their partition function is the Dedekind zeta function. In \cite{Cor}, the first named author has proven that a global function field isomorphism is determined by an isomorphism of dynamical systems for such fields; a result which is reproven here in a more general context.  

Following the philosophy of the Langlands programme, one may replace the consideration of ($n$-dimensional linear representations of) the absolute Galois group $G_{\Kz}$ by that of (modules over) a Hecke algebra (for $\mathrm{GL}_n$). In \cite{CK}, it is shown that for $n \geq 2$, an isomorphism of adelic Hecke algebras for two number fields, that respects the $L^1$-norm, is the same as an isomorphism of their adele rings. The local case is treated in \cite{K}. 

We believe the following questions are interesting to pursue. Firstly, is there a ''Hom-form'' of our main theorem, in which field homomorphisms correspond in a precise way to topological conjugacies of the dynamical systems? Secondly, in the style of Mochizuki's \emph{absolute} version of anabelian geometry (cf.\ \cite{Mo}), one may ask how to reconstruct a number field from its associated dynamical system (or $L$-series), rather than the \emph{relative} result about reconstructing an isomorphism of fields from an isomorphism of their dynamical systems (or $L$-series).

\section{Preliminaries}

In this section, we set notation and introduce the main object of study. 

A monoid is a semigroup with identity element. If $R$ is a ring, we let $R^*$ denote its group of invertible elements.

Given a global field $\Kz$, we use the word \emph{prime} to denote a prime ideal if $\Kz$ is a number field, and to denote an irreducible effective divisor if $\Kz$ is a global function field. Let $\cP_{\Kz}$ be the set of primes of $\Kz$.  If $\mfp \in \cP_{\Kz}$, let $v_{\mfp}$ be the normalized (additive) valuation corresponding to $\mfp$, $\Kz_{\mfp}$ the local field at $\mfp$, and $\cO_\mfp$ its ring of integers. 

Let $\Az_{\Kz,f}$ be the finite adele ring of $\Kz$, $\widehat{\cO}_{\Kz}$ its ring of finite integral adeles and $\Az_{\Kz,f}^*$ the group of ideles (invertible finite adeles), all with their usual topology. Note that in the function field case, we do not single out infinite places, so that in that case, $\Az_{\Kz,f} = \Az_{\Kz}$ is the adele ring, $\widehat{\cO}_{\Kz}$ is the ring of integral adeles and $\Az_{\Kz,f}^* = \Az_{\Kz}^*$ is the full group of ideles. If $\Kz$ is a number field, we denote by $\cO_{\Kz}$ its ring of integers. If $\Kz$ is a global function field, we denote by $q$ the cardinality of the residue field. 

Let $I_{\Kz}$ be the monoid of non-zero integral ideals/effective divisors of our global field $\Kz$, so $I_{\Kz}$ is generated by $\cP_{\Kz}$. We extend the valuation to ideals: if $\mfm \in I_{\Kz}$ and $\mfp \in \cP_{\Kz}$, we define $v_{\mfp}(\mfm) \in \mathbb{Z}_{\geq 0}$ by requiring $\mfm = \prod_{\mfp \in \cP_{\Kz}} \mfp^{v_{\mfp}(\mfm)}$. Let $N$ be the \emph{norm function} on the monoid $I_{\Kz}$: it is the multiplicative function defined on primes $\mfp$ by $N(\mfp):=\# \cO_{\Kz}/\mfp$ if $\Kz$ is a number field and $N(\mfp) = q^{\deg(\mfp)}$ if $\Kz$ is a function field. Given two global fields $\Kz$ and $\Lz$, we call a monoid homomorphism $\varphi: \ I_{\Kz} \to I_{\Lz}$ \emph{norm-preserving} if $N(\varphi(\mfm)) = N(\mfm)$ for all $\mfm \in I_{\Kz}$.

We have a canonical projection $$(\cdot): \ \Az_{\Kz,f}^* \cap \widehat{\cO}_{\Kz} \to I_{\Kz}, \, (x_{\mfp})_{\mfp} \ma \prod_{\mfp} \mfp^{v_{\mfp}(x_{\mfp})}.$$ A \emph{split for $(\cdot)$} is by definition a monoid homomorphism $s_{\Kz}: \ I_{\Kz} \to \Az_{\Kz,f}^* \cap \widehat{\cO}_{\Kz}$ with the property that $(\cdot) \circ s_{\Kz} = \id_{I_{\Kz}}$, and such that for every prime $\mfp$, $s_{\Kz}(\mfp) = (\ldots, 1, \pi_{\mfp}, 1, \ldots)$ for some uniformizer $\pi_{\mfp} \in \Kz_{\mfp}$ with $v_{\mfp}(\pi_{\mfp}) = 1$.

Let $\Gab_{\Kz}$ be the Galois group of a maximal abelian extension of $\Kz$, as profinite topological group. There is an \emph{Artin reciprocity map} $\Az_{\Kz}^* \rightarrow \Gab_{\Kz}$. In the number field case, we embed $\Az_{\Kz,f}^*$ into the group of ideles $\Az_{\Kz}^*$ via $\Az_{\Kz,f}^* \ni x \ma (1,x) \in \Az_{\Kz}^*$, restrict the Artin reciprocity map to $\Az_{\Kz,f}^*$ and call this restriction $\rec_{\Kz}$. In the function field case, $\rec_{\Kz}$ is just the (full) Artin reciprocity map.

\begin{definition}\label{def:X_K} We define the  \emph{topological monoid} $X_{\Kz} \defeq \Gab_{\Kz} \times_{\widehat{\cO}_{\Kz}^*} \widehat{\cO}_{\Kz}$ as the quotient of $\Gab_{\Kz} \times \widehat{\cO}_{\Kz}$ by the subgroup $\menge{(\rec_{\Kz}(u)^{-1},u)}{u \in \widehat{\cO}_{\Kz}^*}$. For $\alpha \in \Gab_{\Kz}$ and $\eta \in \widehat{\cO}_{\Kz}$, we write $[\alpha,\eta]$ for the equivalence class of $(\alpha,\eta)$ in $X_{\Kz}$. Given a split $s_{\Kz}$, we embed $I_{\Kz}$ as a monoid into $X_{\Kz}$ via $$I_{\Kz} \ni \mfm \ma [\rec_{\Kz}(s_{\Kz}(\mfm))^{-1},s_{\Kz}(\mfm)] \in X_{\Kz}.$$ This embedding does not depend on the choice of the split, as this choice is up to units in $\widehat{\cO}_{\Kz}$. We obtain a \emph{dynamical system} $I_{\Kz} \acts X_{\Kz}$, where $\mfm \in I_{\Kz}$ acts on $X_{\Kz}$ via $X_{\Kz} \ni x \ma \mfm \cdot x \in X_{\Kz}$.
\end{definition}

\begin{proposition}
The space $X_{\Kz}$ is Hausdorff. 
\end{proposition}

\bproof
By combining Proposition 1 of Paragraph 4.1 and Proposition 3 of Paragraph 4.2 of \cite{Bour}, we obtain that it suffices to prove that $\widehat{\cO}_{\Kz}^\ast$ is compact and that the action $\widehat{\cO}_{\Kz}^\ast  \acts (\Gab_{\Kz} \times \widehat{\cO}_{\Kz})$ is continuous. The former is clear: for any $\mfp \in \mathcal{P}_{\Kz}$ the unit group $\cO_{\mfp}^\ast$ is compact, hence the product $\widehat{\cO}_{\Kz}^\ast$ is compact as well. To prove the latter it is sufficient that the maps $\widehat{\cO}_{\Kz}^\ast \times \Gab_{\Kz} \to \Gab_{\Kz}$ and $\widehat{\cO}_{\Kz}^\ast \times \widehat{\cO}_{\Kz} \to \widehat{\cO}_{\Kz}$ are continuous. 

The map $\widehat{\cO}_{\Kz}^\ast \times \Gab_{\Kz} \to \Gab_{\Kz}$ is given by $(u, \sigma) \mapsto \rec_{\Kz}(u)^{-1} \sigma$. As the Artin reciprocity map is continuous and $\Gab_{\Kz}$ is a topological group, this map is continuous.

The map $\widehat{\cO}_{\Kz}^\ast \times \widehat{\cO}_{\Kz} \to \widehat{\cO}_{\Kz}$ is given by $(u, a) \mapsto ua$. The topology on $\widehat{\cO}_{\Kz}^\ast$ is finer than the subspace topology obtained from $\widehat{\cO}_{\Kz}$, which causes the inclusion $\widehat{\cO}_{\Kz}^\ast \hookrightarrow \widehat{\cO}_{\Kz}$ to be continuous, which, combined with the fact that $\widehat{\cO}_{\Kz}$ is a topological group, proves that $(u, a) \mapsto ua$ is continuous.

All the necessary conditions are met: we conclude that $X_{\Kz}$ is Hausdorff.
\eproof

\section{The main theorem}

\btheo
\label{THM}
Let $\Kz$ and $\Lz$ be two global fields. The following are equivalent:
\begin{enumerate}
\item[(i)] $I_{\Kz} \acts X_{\Kz}$ and $I_{\Lz} \acts X_{\Lz}$ are \emph{orbit equivalent}, i.e., there exists a homeomorphism $\Phi: \ X_{\Kz} \cong X_{\Lz}$ with $\Phi(I_{\Kz} \cdot x) = I_{\Lz} \cdot \Phi(x)$ for every $x \in X_{\Kz}$.
\item[(ii)] $I_{\Kz} \acts X_{\Kz}$ and $I_{\Lz} \acts X_{\Lz}$ are \emph{conjugate}, i.e., there exists a homeomorphism $\Phi: \ X_{\Kz} \cong X_{\Lz}$ and a monoid isomorphism $\varphi: \ I_{\Kz} \cong I_{\Lz}$ with $\Phi(\mfm \cdot x) = \varphi(\mfm) \cdot \Phi(x)$ for every $\mfm \in I_{\Kz}$ and $x \in X_{\Kz}$.
\item[(iii)] There exists an isomorphism of topological monoids $\Phi: \ X_{\Kz} \cong X_{\Lz}$ which restricts to a monoid isomorphism $I_{\Kz} \cong I_{\Lz}$.
\item[(iv)] There exists a monoid isomorphism $\varphi: \ I_{\Kz} \cong I_{\Lz}$, an isomorphism of topological groups $\psi: \ \Gab_{\Kz} \cong \Gab_{\Lz}$ and splits $s_{\Kz}: \ I_{\Kz} \to \Az_{\Kz,f}^* \cap \widehat{\cO}_{\Kz}$, $s_{\Lz}: \ I_{\Lz} \to \Az_{\Lz,f}^* \cap \widehat{\cO}_{\Lz}$ such that
\bgln
\label{iv1}
  && \psi(\rec_{\Kz}(\cO_{\mfp}^*)) = \rec_{\Lz}(\cO_{\varphi(\mfp)}^*) \ \text{for every prime} \ \mfp \ \text{of} \ \Kz, \\
\label{iv2}
  && \psi(\rec_{\Kz}(s_{\Kz}(\mfm))) = \rec_{\Lz}(s_{\Lz}(\varphi(\mfm))) \ \text{for all} \ \mfm \in I_{\Kz}.
\egln
\end{enumerate}
\etheo

\bremark The main result of \cite{CLMSS} (proven in a purely number theoretical way) implies that the statements in Theorem \ref{THM} are equivalent to 
\begin{enumerate}
\item [(iv)$_{\rm fin}$] \emph{There exists 
\begin{itemize} 
\item a norm-preserving monoid isomorphism $\varphi: \ I_{\Kz} \cong I_{\Lz}$ that restricts to a bijection between the unramified primes of $\Kz$ and $\Lz$, and 
\item an isomorphism of topological groups $\psi: \ \Gab_{\Kz} \cong \Gab_{\Lz}$,
\end{itemize}
such that for every finite abelian extension $\Kz'=\left(\Kz^{\ab}\right)^N/\Kz$ (with $N$ a subgroup in $G_{\Kz}^{\ab}$) and every prime $\mfp$ of $\Kz$ unramified in $\Kz'$, 
we have
$$ \psi \left( \mathrm{Frob}_{\mfp} \right) = \mathrm{Frob}_{\varphi(\mfp)}. $$}
\item[(v)] \emph{There exists an isomorphism of topological groups $\psi: \ \Gab_{\Kz} \cong \Gab_{\Lz}$ such that there is an equality of associated Dirichlet $L$-series: 
\[
  L(\chi,s) = L(\hpsi(\chi),s) \ \text{for all} \ \chi \in \Hab{\Kz}, \ \ \ {\rm where} \ \hpsi \ {\rm is} \ {\rm given} \ {\rm by} \ \hpsi(\chi) = \chi \circ \psi^{-1}.
  \]
}
\item[(vi)] \emph{The fields $\Kz$ and $\Lz$ are isomorphic.}
\end{enumerate}
\eremark

\bremark
We will see that the statements in the theorem are also equivalent to \emph{
\begin{enumerate}
\item[(i$_{\rm N}$)] $I_{\Kz} \acts X_{\Kz}$ and $I_{\Lz} \acts X_{\Lz}$ are \emph{orbit equivalent in a norm-preserving way}, i.e., there exists a homeomorphism $\Phi: \ X_{\Kz} \cong X_{\Lz}$ with $\Phi(I_{\Kz} \cdot x) = I_{\Lz} \cdot \Phi(x)$ for every $x \in X_{\Kz}$. In addition, for every $\mfm \in I_{\Kz}$ and $x \in X_{\Kz}$ there exists $\mfn \in I_{\Lz}$ with $N(\mfm) = N(\mfn)$ and $\Phi(\mfm \cdot x) = \mfn \cdot \Phi(x)$.
\item[(ii$_{\rm N}$)] $I_{\Kz} \acts X_{\Kz}$ and $I_{\Lz} \acts X_{\Lz}$ are \emph{conjugate in a norm-preserving way}, i.e., there exist a homeomorphism $\Phi: \ X_{\Kz} \cong X_{\Lz}$ and a norm-preserving monoid isomorphism $\varphi: \ I_{\Kz} \cong I_{\Kz}$ with $\Phi(\mfm \cdot x) = \varphi(\mfm) \cdot \Phi(x)$ for all $\mfm \in I_{\Kz}$, $x \in X_{\Kz}$.
\item[(iii$_{\rm N}$)] There exists an isomorphism of topological monoids $\Phi: \ X_{\Kz} \cong X_{\Lz}$ which restricts to a norm-preserving monoid isomorphism $I_{\Kz} \cong I_{\Lz}$.
\end{enumerate}}
\eremark

The proof of the theorem proceeds as follows: (i) $\Rarr$ (ii) is the content of Proposition~\ref{i->ii}, and (ii) $\Rarr$ (iii) follows from Lemma~\ref{ii->iii}. We get (iii) $\Rarr$ (iv) by combining Lemma~\ref{PhiG}, Corollary~\ref{PhiO=O} and Lemma~\ref{Phi[pi]}. Proposition~\ref{reci-->monoid} together with Remark~\ref{reci-->monoid_N} and Proposition~\ref{reci-->N} show (iv) $\Rarr$ (iii$_{\rm N}$), and (iii$_{\rm N}$) $\Rarr$ (ii$_{\rm N}$) $\Rarr$ (i$_{\rm N}$) $\Rarr$ (i) is obvious.  

\section{Ideals are dense}

To simplify notation, we drop the index $\Kz$ in this section, i.e., we write $X$ for $X_{\Kz}$, $\Gab$ for $\Gab_{\Kz}$ and so on.
\blemma
\label{IdenseX}
The monoid $\menge{[\rec(s(\mfm))^{-1},s(\mfm)]}{\mfm \in I}$ is dense in $X = \Gab \times_{\widehat{\cO}^*} \widehat{\cO}$.
\elemma
\bproof
For a prime $\mfp$ and a finite abelian extension $L/K$ unramified at $\mfp$ with corresponding projection map $\pi_L: \Gab \twoheadrightarrow \Gal(L/K)$ we have $\pi_L(\rec(s(\mfp))) = \text{Frob}_{\mfp}$. By Chebotarev's Density Theorem (\cite{Neu}, 13.4) for every element $\sigma \in \Gal(L/K)$ there are infinitely many primes $\mfp$ such that $\text{Frob}_{\mfp} = \sigma$. Hence, for any subset $S \subset \cP$ such that $\cP \backslash S$ is finite, the composition $\pi_L \circ \rec \circ s: S \to \Gal(L/K)$ is surjective, thus $\rec(s(S))$ is dense in $\Gab$. 

As a consequence, we show that the set $\bigcup_{\mfm \in I} \menge{[\alpha,s(\mfm)]}{\alpha \in \Gab}$ is in the closure of the image of $I$ in $X$. 

Choose $\alpha \in \Gab$ and $\mfm \in I$ arbitrary and enumerate the primes in $\cP$ by $\mfp_1, \mfp_2, \cdots$. Set $\cP_i = \{\mfp_i, \mfp_{i+1}, \cdots\}$. Moreover, let $H_1 \subset H_2 \subset \cdots$ be an increasing chain of quotient groups of $\Gab$ that gives rise to an identification $\Gab = \plim H_i$. Denote the corresponding quotient map $\Gab \twoheadrightarrow H_i$ by $\pi_i$. As seen before, $\pi_i(\rec(s(\cP_j))) = H_j$ for all $j \in \Nz$. Specifically, for every $i \in \Nz$ we can find $\mfn_i \in \cP_i$ with 
\[
\pi_i(\rec(s(\mfn_i))^{-1}) = \pi_i(\alpha \, \rec(s(\mfm))).
\]
Hence $\mfm_i \defeq \mfm \mfn_i$ satisfy $\lim\limits_{i \to \infty} \rec(s(\mfm_i))^{-1} = \alpha$. By construction, $\lim\limits_{i \to \infty} s(\mfm_i) = s(\mfm)$. 
We are left to show that $\bigcup_{\mfm \in I} \menge{[\alpha,s(\mfm)]}{\alpha \in \Gab}$ is dense in $X$. Let $[\alpha,\eta]$ be an arbitrary element of $X$, where $\eta = (\eta_\mfp)_{\mfp}$. Write $\eta = u \cdot \theta$, where $u \in \widehat{\cO}^\ast$, $\theta = (\theta_\mfp)_{\mfp}\in \widehat{\cO}$ such that $\theta_\mfp$ is either 0 or an integer power of the uniformizer. 

Define
\[
\mfm_j := \prod_{i = 1}^j \mfp_i^{\min(v_{\mfp}(\theta_\mfp), i)},
\]
where $\min(\infty, i) = i$ for all $i$. Then $s(\mfm_j)_\mfp$ converges to $\theta_\mfp$ for all $\mfp$, hence $\lim\limits_{j \to \infty} s(\mfm_j) = \theta$. It follows that
$$\lim\limits_{j \to \infty} [\alpha \, \rec(u),s(\mfm_j)] =  [\alpha \, \rec(u), \theta] =[\alpha \, \rec(u),u^{-1} \eta] = [\alpha,\eta].$$
This concludes the proof.
\eproof

\section{From orbit equivalence to conjugacy}
Let $\Kz$ and $\Lz$ be two global fields.
\bprop
\label{i->ii}
The following are equivalent:
\begin{itemize}
\item[(a)] $I_{\Kz} \acts X_{\Kz}$ and $I_{\Lz} \acts X_{\Lz}$ are orbit equivalent. 
\item[(b)] $I_{\Kz} \acts X_{\Kz}$ and $I_{\Lz} \acts X_{\Lz}$ are conjugate. 
\end{itemize}
If, in addition to (a), for every $\mfm \in I_{\Kz}$ and $x \in X_{\Kz}$, we can always choose $\mfn \in I_{\Lz}$ so that $N(\mfn) = N(\mfm)$ and $\Phi(\mfm \cdot x) = \mfn \cdot \Phi(x)$, then we can also find a norm-preserving monoid isomorphism in (b).
\eprop
\bproof
It obviously suffices to show \an{(a) $\Rarr$ (b)}. First of all, let $1 = [1,1] \in X_{\Kz}$ be the unit of $X_{\Kz}$. We claim that $\Phi(1)$ is invertible in the monoid $X_{\Lz}$. Indeed, since $I_{\Kz} = I_{\Kz} \cdot 1$ is dense in $X_{\Kz}$, we know that $I_{\Lz} \cdot \Phi(1) = \Phi(I_{\Kz} \cdot 1)$ is dense in $X_{\Lz}$. Now assume that $\Phi(1)$ is not invertible in $X_{\Lz}$. This means that $\Phi(1) = [\alpha,\eta]$, and there exists a prime $\mfq$ of $\Lz$ with $v_{\mfq}(\eta_{\mfq}) \geq 1$. But then, every $[\beta,\xi] \in \overline{I_{\Lz} \cdot \Phi(1)}$ satisfies $v_{\mfq}(\xi_{\mfq}) \geq 1$. This contradicts $\overline{I_{\Lz} \cdot \Phi(1)} = X_{\Lz}$. So $\Phi(1)$ must be invertible in $X_{\Lz}$.

By (a), we know that for every $\mfm \in I_{\Kz}$, there exists $\mfn \in I_{\Lz}$ with $\Phi(\mfm \cdot 1) = \mfn \cdot \Phi(1)$. Since $\Phi(1)$ is invertible, $\mfn$ is uniquely determined. In other words, if $\mfn_1, \mfn_2 \in I_L$ both satisfy $\Phi(\mfm \cdot 1) = \mfn_1 \cdot \Phi(1)$ and $\Phi(\mfm \cdot 1) = \mfn_2 \cdot \Phi(1)$, then this implies $\mfn_1 \cdot \Phi(1) = \mfn_2 \cdot \Phi(1)$ and hence $\mfn_1 = \mfn_2$. Let $\varphi: \ I_{\Kz} \to I_{\Lz}$ be the map with $\Phi(\mfm \cdot 1) = \varphi(\mfm) \cdot \Phi(1)$ for every $\mfm \in I_{\Kz}$. Similarly, let $\phi: \ I_{\Lz} \to I_{\Kz}$ be the map uniquely determined by $\Phi^{-1}(\mfn \cdot \Phi(1)) = \phi(\mfn)$ for all $\mfn \in I_{\Lz}$. It is obvious that $$\phi \circ \varphi = \id_{I_{\Kz}} \mbox{ and }\varphi \circ \phi = \id_{I_{\Lz}}.$$ So $\phi = \varphi^{-1}$, and in particular, $\varphi$ is a bijection.

Our first step is to show that for every prime $\mfp$ of $\Kz$ and $\mfm \in I_{\Kz}$, there exists a prime $\mfq$ of $\Lz$ such that $\varphi(\mfp \cdot \mfm) = \mfq \cdot \varphi(\mfm)$. At this point, it might be that $\mfq$ depends on $\mfm$, but we will see later on that this is not the case, i.e., that $\mfq = \varphi(\mfp)$. Looking at $$\varphi(\mfp \cdot \mfm) \cdot \Phi(1) = \Phi(\mfp \cdot \mfm) \in I_{\Lz} \cdot \Phi(\mfm) = I_{\Lz} \cdot \varphi(\mfm) \cdot \Phi(1),$$ we conclude $\varphi(\mfp \cdot \mfm) \in I_{\Lz} \cdot \varphi(\mfm)$ since $\Phi(1)$ is invertible. So there exists $\mfa \in I_{\Lz}$ with $\varphi(\mfp \cdot \mfm) = \mfa \cdot \varphi(\mfm)$. If $\mfa$ is not prime, then there exist $\mfb, \mfc \in I_{\Lz}$ with $\mfb \neq \cO_{\Lz}$, $\mfc \neq \cO_{\Lz}$ and $\mfa = \mfb \cdot \mfc$. Because of orbit equivalence, we must have
\bglnoz
  \mfp \cdot \mfm &=& \Phi^{-1}(\varphi(\mfp \cdot \mfm) \cdot \Phi(1)) = \Phi^{-1}(\mfa \cdot \varphi(\mfm) \cdot \Phi(1)) 
  = \Phi^{-1}(\mfb \cdot \mfc \cdot \varphi(\mfm) \cdot \Phi(1)) \\
  &=& \ti{\mfb} \cdot \Phi^{-1}(\mfc \cdot \varphi(\mfm) \cdot \Phi(1)) = \ti{\mfb} \cdot \ti{\mfc} \cdot \Phi^{-1}(\varphi(\mfm) \cdot \Phi(1))
  = \ti{\mfb} \cdot \ti{\mfc} \cdot \mfm
\eglnoz
for some $\ti{\mfb}, \ti{\mfc} \in I_{\Kz}$. Let $$x = \Phi^{-1}(\mfc \cdot \varphi(\mfm) \cdot \Phi(1)) \mbox{ and } y = \Phi(x) = \mfc \cdot \varphi(\mfm) \cdot \Phi(1).$$ Let $y = [\beta,\xi]$. As $\mfb \neq \cO_{\Lz}$, there exists a prime $\mfq$ of $\Lz$ with $v_{\mfq}(\mfb) \geq 1$. Obviously, $v_{\mfq}(\xi_{\mfq}) \neq \infty$, and we conclude that $v_{\mfq}(s_{\Lz}(\mfb)_{\mfq} \cdot \xi_{\mfq}) \neq v_{\mfq}(\xi_{\mfq})$, hence $\mfb \cdot y \neq y$. Therefore, $x \neq \ti{\mfb} \cdot x$, and thus $\ti{\mfb} \neq \cO_{\Kz}$. Similarly, $\ti{\mfc} \neq \cO_{\Kz}$. But then, it is impossible to have $\mfp \cdot \mfm = \ti{\mfb} \cdot \ti{\mfc} \cdot \mfm$ in $I_{\Kz}$. Hence, $\mfa$ must be prime. In particular, $\varphi$ induces a bijective map on the set of primes.

Our next step is to show $\varphi(\mfp^n) = \varphi(\mfp)^n$ for every prime $\mfp$ of $\Kz$ and $n \in \Nz$. By our previous step, there exist primes $\mfq_1, \dotsc, \mfq_n$ of $\Lz$ such that $\varphi(\mfp^n) = \mfq_1 \cdots \mfq_n$. Thus $$I_{\Lz} \cdot \varphi(\mfp^n) \cdot \Phi(1) \subseteq I_{\Lz} \cdot \mfq_i \cdot \Phi(1)$$ for every $1 \leq i \leq n$. Applying $\Phi^{-1}$, this implies that $I_{\Kz} \cdot \mfp^n \subseteq I_{\Kz} \varphi^{-1}(\mfq_i)$. But the only prime $\ti{\mfp}$ of $\Kz$ with $I_{\Kz} \cdot \mfp^n \subseteq I_{\Kz} \cdot \ti{\mfp}$ is $\mfp$. We therefore conclude that $\mfq_i = \varphi(\mfp)$ for all $1 \leq i \leq n$.

Let us now show that for every $\mfm \in I_{\Kz}$, we have $v_{\varphi(\mfp)}(\varphi(\mfm)) = v_{\mfp}(\mfm)$. It is clear that $v_{\varphi(\mfp)}(\varphi(\mfm))$ is the unique $v \in \Nz$ with the properties $$I_{\Lz} \cdot \varphi(\mfm) \subseteq I_{\Lz} \cdot \varphi(\mfp)^v = I_{\Lz} \cdot \varphi(\mfp^v) \mbox{ and }I_{\Lz} \cdot \varphi(\mfm) \nsubseteq I_{\Lz} \cdot \varphi(\mfp)^{v+1} = I_{\Lz} \cdot \varphi(\mfp^{v+1}).$$ Here we used the previous step. Multiplying the equations with $\Phi(1)$ and applying $\Phi^{-1}$, these two properties transform into $$I_{\Kz} \cdot \mfm \subseteq I_{\Kz} \cdot \mfp^v \mbox{ and }I_{\Kz} \cdot \mfm \nsubseteq I_{\Kz} \cdot \mfp^{v+1}.$$ Thus $v = v_{\mfp}(\mfm)$. In particular, this shows that $\varphi$ is a monoid isomorphism.

We can now finish the proof for \an{(a) $\Rarr$ (b)}: For every $\mfm, \mfn \in I_{\Kz}$, we have $$\Phi(\mfm \cdot \mfn) = \varphi(\mfm \cdot \mfn) \cdot \Phi(1) = \varphi(\mfm) \cdot \varphi(\mfn) \cdot \Phi(1) = \varphi(\mfm) \cdot \Phi(\mfn).$$ As $I_{\Kz}$ is dense in $X_{\Kz}$, we conclude that $\Phi(\mfm \cdot x) = \varphi(\mfm) \cdot \Phi(x)$ for every $\mfm \in I_{\Kz}$ and $x \in X_{\Kz}$.

Let us now assume that, in addition to (a), we can, for every $\mfm \in I_{\Kz}$ and $x \in X_{\Kz}$, choose $\mfn \in I_{\Lz}$ so that $N(\mfn) = N(\mfm)$ and $\Phi(\mfm \cdot x) = \mfn \cdot \Phi(x)$. Taking $x=1$, we see that $N(\varphi(\mfm)) = N(\mfm)$.
\eproof

\section{From conjugacy to isomorphism of topological monoids}

\blemma
\label{ii->iii}
Assume that $I_{\Kz} \acts X_{\Kz}$ and $I_{\Lz} \acts X_{\Lz}$ are conjugate.  
Then $\Phi(1)$ is invertible in $X_{\Lz}$ and $\ti{\Phi}: \ X_{\Kz} \to X_{\Lz}, \, x \ma \Phi(x) \cdot \Phi(1)^{-1}$ is an isomorphism of topological monoids satisfying $\ti{\Phi}(\mfm) = \varphi(\mfm)$ for every $\mfm \in I_{\Kz}$.
\elemma
\bproof
In the proof of the previous proposition, we have already seen that $\Phi(1)$ is invertible in $X_{\Lz}$. Obviously, $\ti{\Phi}$ defined by $\ti{\Phi}(x) = \Phi(x) \cdot \Phi(1)^{-1}$ is still a homeomorphism. Moreover, for every $\mfm, \mfn \in I_K$, we have \begin{align*} \ti{\Phi}(\mfm \cdot \mfn) &= \Phi(\mfm \cdot \mfn) \Phi(1)^{-1} \\ &= \varphi(\mfm) \cdot \varphi(\mfn) \cdot \Phi(1) \cdot \Phi(1)^{-1} \\ &= (\varphi(\mfm) \cdot \Phi(1) \cdot \Phi(1)^{-1}) \cdot (\varphi(\mfn) \cdot \Phi(1) \cdot \Phi(1)^{-1}) \\ &= \ti{\Phi}(\mfm) \cdot \ti{\Phi}(\mfn).\end{align*} As $I_{\Kz}$ is dense in $X_{\Kz}$, we conclude that $\ti{\Phi}$ is a monoid homomorphism. Finally, we have for every $\mfm \in I_{\Kz}$: $$\ti{\Phi}(\mfm) = \Phi(\mfm \cdot 1) \cdot \Phi(1)^{-1} = \varphi(\mfm) \cdot \Phi(1) \cdot \Phi(1)^{-1} = \varphi(\mfm).$$
\eproof

\bremark
Note that in Lemma~\ref{ii->iii}, the monoid isomorphism $\varphi$ stays the same. In particular, if it was norm-preserving, then we will keep that property.
\eremark

\section{From conjugacy to reciprocity isomorphism}

\blemma
\label{PhiG}
Let $\Phi$ be an isomorphism of topological monoids $\Phi: \ X_{\Kz} \cong X_{\Lz}$. Then $\Phi$ restricts to an isomorphism of topological groups $\Gab_{\Kz} \cong \Gab_{\Lz}$.
\elemma
\bproof
This is clear as $\Gab_{\Kz}$ is the group of invertible elements in $X_{\Kz}$.
\eproof

\begin{notation} Let us denote the restriction of $\Phi$ to $\Gab_{\Kz}$ again by $\Phi$. For $\mfm \in I_{\Kz}$, let $\cP(\mfm) = \menge{\mfp \in \cP_{\Kz}}{\mfp | \mfm}$. Given $S \subseteq \cP_{\Kz}$, let $1_S$ be the element of $\widehat{\cO}_{\Kz}$ with $(1_S)_{\mfp} = 1$ if $\mfp \in S$ and $(1_S)_{\mfp} = 0$ if $\mfp \notin S$. Also, let $e$ be the identity in $\Gab_{\Kz}$. We do not distinguish in our notation between $e \in \Gab_{\Kz}$ and $e \in \Gab_{\Lz}$ since the meaning becomes clear from the context.
\end{notation}
\blemma
\label{Phi1_S}
Assume that $\Phi$ is an isomorphism of topological monoids $\Phi: \ X_{\Kz} \cong X_{\Lz}$ which restricts to an isomorphism $\varphi: \ I_{\Kz} \cong I_{\Lz}$. We denote the bijection $\cP_{\Kz} \cong \cP_{\Lz}$ that $\varphi$ induces also by $\varphi$. Then, for every $S \subseteq \cP_{\Kz}$, we have $\Phi([e,1_S]) = [e,1_{\varphi(S)}]$.
\elemma
\bproof
Applying $\Phi$ to both sides of the equation
$$[e,1_S] \cdot X_{\Kz} = \bigcap_{\mfm \in I_{\Kz}, \, \cP(\mfm) \cap S = \emptyset} \mfm \cdot X_{\Kz},$$
we obtain
$$\Phi([e,1_S]) \cdot X_{\Lz} = \bigcap_{\mfm \in I_{\Kz}, \, \cP(\mfm) \cap S = \emptyset} \varphi(\mfm) \cdot X_{\Lz} = [e,1_{\varphi(S)}] \cdot X_{\Lz}.$$
It follows that if $\Phi([e,1_S]) = [\beta,\xi]$, then $\xi_{\mfq} = 0$ if $\mfq \notin \varphi(S)$, and $\xi_{\mfq} \in {\widehat{\cO}_{\mfq}}^*$ if $\mfq \in \varphi(S)$. As $[e,1_S]^2 = [e,1_S]$ and since $\Phi$ is a monoid homomorphism, we must have $$\Phi([e,1_S])^2 = \Phi([e,1_S]).$$ Thus there exists $u \in {\widehat{\cO}_{\Lz}}^*$ with $\beta^2 = \beta \, \rec_{\Lz}(u)^{-1}$ and $\xi^2 = u \xi$. We deduce that $\beta = \rec_{\Lz}(u)^{-1}$ and $\xi = u \cdot 1_{\varphi(S)}$. Hence
$$
  \Phi([e,1_S]) = [\beta,\xi] = [\rec_{\Lz}(u)^{-1},u \cdot 1_{\varphi(S)}] = [e,1_{\varphi(S)}].
$$
\eproof
\begin{notation} With the same notation as above, define for $S \subseteq \cP_{\Kz}$ the subgroup
$$
N_S = \rec_{\Kz}(\prod_{\mfp \in S} \cO_{\mfp}^*) 
\subseteq \Gab_{\Kz}.
$$
We view $\prod_{\mfp \in S} \cO_{\mfp}^*$ as a subgroup of $\widehat{\cO}_{\Kz}^*$ via the embedding $v \ma (\dotsc, 1, v, 1, \dotsc)$.  
\end{notation}
\blemma
\label{PhiN}
Let $\Phi$ be an isomorphism of topological monoids $\Phi: \ X_{\Kz} \cong X_{\Lz}$ which restricts to a isomorphism $\varphi: \ I_{\Kz} \cong I_{\Lz}$. Then $\Phi(N_S) = N_{\varphi(S)}$.
\elemma
\bproof
Let $S^c = \cP_{\Kz} \setminus S$ and write $\mu_{S^c}$ for the monoid homomorphism $$ \mu_{S^c} \colon \Gab_{\Kz} \to X_{\Kz} \cdot [e,1_{S^c}], \, \alpha \ma \alpha \cdot [e,1_{S^c}].$$ Define $\mu_{\varphi(S)^c}$ in an analogous way. Because of the previous two lemmas, $\Phi$ restricts to isomorphisms $\Gab_{\Kz} \cong \Gab_{\Lz}$ and $X_{\Kz} \cdot [e,1_{S^c}] \cong X_{\Lz} \cdot [e,1_{\varphi(S)^c}]$ such that the diagram
\bgloz
  \xymatrix@C=24mm{
  \Gab_{\Kz} \ar[d]^{\cong}_{\Phi} \ar[r]^{\mu_{S^c}} & X_{\Kz} \cdot [e,1_{S^c}] \ar[d]^{\Phi}_{\cong} \\ 
  \Gab_{\Lz} \ar[r]^{\mu_{\varphi(S)^c}} & X_{\Lz} \cdot [e,1_{\varphi(S)^c}]
  }
\egloz
commutes. Hence it follows that $$\Phi((\mu_{S^c})^{-1}([e,1_{S^c}])) = (\mu_{\varphi(S)^c})^{-1}([e,1_{\varphi(S)^c}]).$$ Now $[\alpha,1] \in \Gab_{\Kz}$ lies in $(\mu_{S^c})^{-1}([e,1_{S^c}])$ if and only if $[\alpha,1_{S^c}] = [e,1_{S^c}]$ if and only if there exists $u \in \widehat{\cO}_{\Kz}^*$ with $\alpha = \rec_{\Kz}(u)^{-1}$ and $u \cdot 1_{S^c} = 1_{S^c}$. The latter equation means that $u_{\mfp} = 1$ if $\mfp \in S^c$, or in other words, $u \in \prod_{\mfp \in S} \cO_{\mfp}^*$. Thus, $$(\mu_{S^c})^{-1}([e,1_{S^c}]) = N_S.$$ In an entirely analogous way, we see that $(\mu_{\varphi(S)^c})^{-1}([e,1_{\varphi(S)^c}]) = N_{\varphi(S)}$.
\eproof
\bcor
\label{PhiO=O}
We have an identification $\Phi(\rec_{\Kz}( \cO_{\mfp}^*)) = \rec_{\Lz}( \cO_{\varphi(\mfp)}^*)$.\qed
\ecor

\begin{notation} For $\mfp \in \cP_{\Kz}$, let $\pi_{\mfp} \in \widehat{\cO}_{\Kz}$ be such that $(\pi_{\mfp})_{\mfq} = 1$ if $\mfq \neq \mfp$ and $v_{\mfp}(\pi_{\mfp}) = 1$.\end{notation}
\blemma
\label{Phipi}
Assume that $\Phi$ is an isomorphism of topological monoids $\Phi: \ X_{\Kz} \cong X_{\Lz}$ which restricts to a isomorphism $\varphi: \ I_{\Kz} \cong I_{\Lz}$. Then for every $\mfp \in \cP_{\Kz}$ there exists an element $\pi_{\varphi(\mfp)} \in \widehat{\cO}_{\Lz}$ with $$ \left\{ \begin{array}{l} (\pi_{\varphi(\mfp)})_{\mfr} = 1 \mbox{ if }\mfr \neq \varphi(\mfp), \\v_{\varphi(\mfp)}(\pi_{\varphi(\mfp)}) = 1, \\ \Phi([e,\pi_{\mfp}]) = [e,\pi_{\varphi(\mfp)}]. \end{array} \right. $$ 
\elemma
\bproof
We have
$$ \Phi(\menge{[\alpha,\eta]}{v_{\mfp}(\eta_{\mfp}) > 0}) = \menge{[\beta,\xi]}{v_{\varphi(\mfp)}(\xi_{\varphi(\mfp)}) > 0}. $$
The reason is that for $[\alpha,\eta] \in X_{\Kz}$, $$v_{\mfp}(\eta_{\mfp}) = 0 \LRarr \eta_{\mfp} \in \cO_{\mfp}^* \LRarr [\alpha,\eta] \cdot [e,1_{\gekl{\mfp}}] \cdot X_{\Kz} = [e,1_{\gekl{\mfp}}] \cdot X_{\Kz}.$$ The latter condition is equivalent to $\Phi([\alpha,\eta]) \cdot [e,1_{\gekl{\varphi(\mfp)}}] \cdot X_{\Lz} = [e,1_{\gekl{\varphi(\mfp)}}] \cdot X_{\Lz}$, and this in turn means that we can write $\Phi([\alpha,\eta]) = [\beta,\xi]$ for some $\xi \in \widehat{\cO}_{\Lz}$ with $v_{\varphi(\mfp)}(\xi_{\varphi(\mfp)}) = 0$. Here we used Lemma~\ref{Phi1_S}.

As we obviously have $[e,\pi_{\mfp}] \cdot X_{\Kz} = \menge{[\alpha,\eta]}{v_{\mfp}(\eta_{\mfp}) > 0}$, we conclude that $$\Phi([e,\pi_{\mfp}]) \cdot X_{\Lz} = \menge{[\beta,\xi]}{v_{\varphi(\mfp)}(\xi_{\varphi(\mfp)}) > 0}.$$ So we can write $\Phi([e,\pi_{\mfp}]) = [\beta,\xi]$ for some $\xi \in \widehat{\cO}_{\Lz}$ with $v_{\varphi(\mfp)}(\xi_{\varphi(\mfp)}) = 1$.

Moreover, we have $[e,\pi_{\mfp}] \cdot [e,1_{\gekl{\mfp}^c}] = [e,1_{\gekl{\mfp}^c}]$. Hence, by Lemma~\ref{Phi1_S}, we conclude that $$[\beta,\xi \cdot 1_{\gekl{\varphi(\mfp)}^c}] = [\beta,\xi] \cdot [e,1_{\gekl{\varphi(\mfp)}^c}] = [e,1_{\gekl{\varphi(\mfp)}^c}].$$ Therefore, there exists $u \in \widehat{\cO}_{\Lz}^*$ with $\beta = \rec_{\Lz}(u)^{-1}$ and $\xi \cdot 1_{\gekl{\varphi(\mfp)}^c} = u \cdot 1_{\gekl{\varphi(\mfp)}^c}$. The latter equation implies that $\xi_{\mfr} = u_{\mfr} \in \widehat{\cO}_{\mfr}^*$ if $\mfr \neq \varphi(\mfp)$. As we have shown above that $v_{\varphi(\mfp)}(\xi_{\varphi(\mfp)}) = 1$, we conclude that $\pi_{\varphi(\mfp)} = u^{-1} \xi$ satisfies $(\pi_{\varphi(\mfp)})_{\mfr} = 1$ if $\mfr \neq \varphi(\mfp)$ and $v_{\varphi(\mfp)}(\pi_{\varphi(\mfp)}) = 1$. Also, $$\Phi([e,\pi_{\mfp}]) = [\beta,\xi] = [\rec_{\Lz}(u)^{-1},u \pi_{\varphi(\mfp)}] = [e,\pi_{\varphi(\mfp)}],$$ as desired.
\eproof

\blemma
\label{Phi[pi]}
Assume that $\Phi$ is an isomorphism of topological monoids $\Phi: \ X_{\Kz} \cong X_{\Lz}$ which restricts to a isomorphism $\varphi: \ I_{\Kz} \cong I_{\Lz}$. Given $\mfp \in \cP_{\Kz}$, let $\pi_{\mfp}$ be defined as above, and choose $\pi_{\varphi(\mfp)} \in \widehat{\cO}_{\Lz}$ as in Lemma~\ref{Phipi}. Then $\Phi(\rec_{\Kz}(\pi_{\mfp})) = \rec_{\Lz}(\pi_{\varphi(\mfp)})$.
\elemma
\bproof
We know that $\Phi$ sends $\mfp \in X_{\Kz}$ to $\varphi(\mfp) \in X_{\Lz}$. At the same time, $$\mfp = [\rec_{\Kz}(\pi_{\mfp})^{-1},\pi_{\mfp}] = [\rec_{\Kz}(\pi_{\mfp})^{-1},1] \cdot [e,\pi_{\mfp}].$$ Let $\beta \in \Gab_{\Lz}$ satisfy $\Phi([\rec_{\Kz}(\pi_{\mfp})^{-1},1]) = [\beta,1]$. Then, by Lemma~\ref{Phipi}, $\Phi(\mfp) = [\beta,\pi_{\varphi(\mfp)}]$. Hence $\varphi(\mfp) = [\rec_{\Lz}(\pi_{\varphi(\mfp)})^{-1},\pi_{\varphi(\mfp)}]$ is equal to $[\beta,\pi_{\varphi(\mfp)}]$. Therefore, there exists $u \in \widehat{\cO}_{\Lz}^*$ with $\rec_{\Lz}(\pi_{\varphi(\mfp)})^{-1} = \beta \rec_{\Lz}(u)^{-1}$ and $\pi_{\varphi(\mfp)} = u \pi_{\varphi(\mfp)}$. But the latter equality implies $u=1$. Hence $\beta = \rec_{\Lz}(\pi_{\varphi(\mfp)})^{-1}$. We conclude that $$\Phi(\rec_{\Kz}(\pi_{\mfp})^{-1}) = \Phi([\rec_{\Kz}(\pi_{\mfp})^{-1},1]) = [\beta,1] = [\rec_{\Lz}(\pi_{\varphi(\mfp)})^{-1},1] = \rec_{\Lz}(\pi_{\varphi(\mfp)})^{-1},$$ and thus $\Phi(\rec_{\Kz}(\pi_{\mfp})) = \rec_{\Lz}(\pi_{\varphi(\mfp)})$.
\eproof

\section{From reciprocity isomorphism to isomorphism of topological monoids}
\label{Reciprocity-->Monoids}

In this section, we start by assuming condition (iv) from the main theorem \ref{THM}. 

\bprop
\label{reci-->monoid}
Assume \textup{(iv)}. 
Then there exists an isomorphism of topological monoids $\Phi: \ X_{\Kz} \cong X_{\Lz}$ which restricts to the monoid isomorphism $\varphi: \ I_{\Kz} \cong I_{\Lz}$.
\eprop
\bproof
For $\mfm \in I_{\Kz}$, define the embedding $$\iota_{\mfm}: \ \Gab_{\Kz} \to X_{\Kz}, \, \alpha \ma [\alpha,s_K(\mfm)].$$ For $\mfm_1 \neq \mfm_2$, it is obvious that $\iota_{\mfm_1}$ and $\iota_{\mfm_2}$ have disjoint images. Define $\Phi$ on $\iota_{\mfm}(\Gab_{\Kz})$ by $$\Phi([\alpha,s_{\Kz}(\mfm)]) \defeq [\psi(\alpha),s_{\Lz}(\varphi(\mfm))].$$ Let us show that $\Phi$ has a unique continuous extension to $X_{\Kz}$. Uniqueness is clear as $\bigcup_{\mfm \in I_{\Kz}} \iota_{\mfm}(\Gab_{\Kz})$ is dense in $X_{\Kz}$ (see \S~\ref{IdenseX}).

For $[\alpha,\eta] \in X_{\Kz}$, let $(\alpha_i)_i$ and $(\mfm_i)_i$ be sequences in $\Gab_{\Kz}$ and $I_{\Kz}$, respectively, such that $\lim\limits_{i \to \infty} [\alpha_i,s_K(\mfm_i)] = [\alpha,\eta]$. Given $\mfp \in \cP_{\Kz}$, either $\eta_{\mfp} = 0$, in which case $$\lim\limits_{i \to \infty} v_{\mfp}(\mfm_i) = \lim\limits_{i \to \infty} v_{\mfp}(s_{\Kz}(\mfm_i)_{\mfp}) = \infty,$$ or $\eta_{\mfp} \neq 0$, which implies that $v_{\mfp}(\mfm_i) = v_{\mfp}(s_{\Kz}(\mfm_i)_{\mfp}) = v_{\mfp}(\eta_{\mfp})$ for sufficiently large $i$. The latter equality means that $s_{\Kz}(\mfm_i)_{\mfp} = s_{\Kz}(\mfp^{v_{\mfp}(\eta_{\mfp})})_{\mfp}$ for sufficiently large $i$, or in other words, $s_{\Kz}(\mfm_i)_{\mfp}$ is eventually constant. Therefore, if we define $\eta' \in \widehat{\cO}_{\Kz}$ by setting $$ \eta'_{\mfp} = \left\{ \begin{array}{l} 0 \mbox{ if }\eta_{\mfp} = 0; \\ s_{\Kz}(\mfp^{v_{\mfp}(\eta_{\mfp})})_{\mfp} \mbox{ if } \eta_{\mfp} \neq 0, \end{array} \right. $$ then $\lim\limits_{i \to \infty} s_{\Kz}(\mfm_i) = \eta'$. Thus, for $\mfq \in \cP_L$,
\begin{itemize} \item either $\eta_{\varphi^{-1}(\mfq)} = 0$, in which case $\lim\limits_{i \to \infty} v_{\mfq}(s_{\Lz}(\varphi(\mfm_i))_{\mfq}) = \infty$,
\item or $\eta_{\varphi^{-1}(\mfq)} \neq 0$, in which case $v_{\mfq}(s_{\Lz}(\varphi(\mfm_i))_{\mfq})$ is eventually constant, namely equal to $v_{\varphi^{-1}(\mfq)}(\eta_{\varphi^{-1}(\mfq)})$. 
\end{itemize}
Therefore, defining $\xi \in \widehat{\cO}_{\Lz}$ by setting $$\xi_{\mfq} = \left\{ \begin{array}{l} 0, \mbox{ if }\eta_{\varphi^{-1}(\mfq)} = 0, \\ s_{\Lz}(\mfq^\lambda)_{\mfq} \mbox{ with } \lambda:={v_{\varphi^{-1}(\mfq)}(\eta_{\varphi^{-1}(\mfq)})}, \mbox{ if }\eta_{\varphi^{-1}(\mfq)} \neq 0, \end{array} \right. $$ we see that $\lim\limits_{i \to \infty} s_{\Lz}(\varphi(\mfm_i)) = \xi$. 

Set $S \defeq \menge{\mfp \in \cP_K}{\eta_{\mfp} = 0}$ and define $v = (v_{\mfp})_{\mfp} \in \prod_{\mfp \in S^c} \cO_{\mfp}^*$ by the equation $\eta'_{\mfp} = v_{\mfp} \cdot \eta_{\mfp}$ for $\mfp \in S^c$.

Let $N$ be the subgroup $N:=\rec_{\Kz}(\prod_{\mfp \in S} \cO_{\mfp}^*)$ of $\Gab_{\Kz}$. Let us show that $(\alpha_i)_i$ converges in $\Gab_{\Kz} / N$. Let $(\alpha_{i_k})_k$ be a convergent subsequence of $(\alpha_i)_i$ with $\lim\limits_{k \to \infty} \alpha_{i_k} = \alpha'$. Then $$\lim\limits_{k \to \infty} [\alpha_{i_k},s_{\Kz}(\mfm_{i_k})] = [\alpha',\eta'] = [\alpha,\eta].$$ So there exists $u \in \widehat{\cO}_{\Kz}^*$ with $\alpha' = \alpha \cdot \rec_{\Kz}(u)^{-1}$ and $\eta' = u \eta$. The latter equation implies $u_{\mfp} = v_{\mfp}$ for all $\mfp \in S^c$, and thus $\alpha' = \alpha \cdot \rec_{\Kz}(v)^{-1}$ in $\Gab_{\Kz} / N$. As our convergent subsequence $(\alpha_{i_k})_k$ was arbitrary, we indeed have $$\lim\limits_{i \to \infty} \alpha_i = \alpha \cdot \rec_{\Kz}(v)^{-1} \mbox{ in } \Gab_{\Kz} / N.$$ Hence $$\lim\limits_{i \to \infty} \psi(\alpha_i) = \psi(\alpha) \cdot \psi(\rec_{\Kz}(v)^{-1}) \mbox{ in }\Gab_{\Lz} / \psi(N).$$ Thus $$\lim\limits_{i \to \infty} [\psi(\alpha_i),s_{\Lz}(\varphi(\mfm_i))] = [\psi(\alpha) \cdot \psi(\rec_{\Kz}(v)^{-1}),\xi].$$ Namely, if $(\psi(\alpha_{i_k}))_k$ is a convergent subsequence of $(\psi(\alpha_i))_i$, then $$\lim\limits_{k \to \infty} \psi(\alpha_{i_k}) = \psi(\alpha) \cdot \psi(\rec_{\Kz}(v)^{-1})\cdot  \rec_{\Lz}(u) \mbox{ in }\Gab_{\Lz},$$ for some $u \in \prod_{\mfq \in \varphi(S)} {\cO}_{\mfq}^*$, by \eqref{iv1}. So
\bglnoz
  \lim\limits_{i \to \infty} [\psi(\alpha_i),s_{\Lz}(\varphi(\mfm_i))]
  &=& [\psi(\alpha) \cdot \psi(\rec_{\Kz}(v)^{-1}) \cdot \rec_{\Lz}(u),\xi] \\
  &=& [\psi(\alpha) \cdot \psi(\rec_{\Kz}(v)^{-1}),u \xi]\\
  &=& [\psi(\alpha) \cdot \psi(\rec_{\Kz}(v)^{-1}),\xi].
\eglnoz
This shows that $\Phi$ extends continuously to a map $\Gab_{\Kz} \to X_{\Kz}$, again denoted by $\Phi$.

Let us show that $\Phi$ is continuous. Suppose that $\lim\limits_{i \to \infty} [\alpha_i,\eta_i] = [\alpha,\eta]$ in $X_{\Kz}$. Define $\eta'$, $\xi$ and $v$ for $\eta$ as before, and define $\eta'_i$, $\xi_i$ and $v_i$ for $\eta_i$. We then have $$\Phi([\alpha_i,\eta_i]) = [\psi(\alpha_i) \cdot \psi(\rec_{\Kz}(v_i)^{-1}),\xi_i].$$ We want to show that $$\lim\limits_{i \to \infty} [\psi(\alpha_i) \cdot \psi(\rec_{\Kz}(v_i)^{-1}),\xi_i] = \Phi([\alpha,\eta]).$$ It is clear that $\lim\limits_{i \to \infty} \eta'_i = \eta'$, and also $\lim\limits_{i \to \infty} \xi_i = \xi$. Passing to convergent subsequences if necessary, we may assume that $$\lim\limits_{i \to \infty} \eta_i = \ti{\eta}; \ \  \lim\limits_{i \to \infty} v_i = v \mbox{ and }\lim\limits_{i \to \infty} \alpha_i = \ti{\alpha}.$$ Obviously, there exists $u \in \widehat{\cO}^*_{\Kz}$ with $\ti{\eta} = u \eta$. Therefore, $$\eta' = \lim\limits_{i \to \infty} \eta'_i = \lim\limits_{i \to \infty} v_i \eta_i = v u \eta.$$ As we have seen before, we must have $$\Phi([\alpha,\eta]) = [\psi(\alpha) \cdot \psi(\rec_{\Kz}(vu)^{-1}),\xi].$$ Now let $S \defeq \menge{\mfp \in \cP_{\Kz}}{\eta_{\mfp} = 0},$ and let $N = \rec_{\Kz}(\prod_{\mfp \in S^c} \cO_{\mfp}^*)$. Obviously, $u$ is uniquely determined modulo $\prod_{\mfp \in S^c} \cO_{\mfp}^*$. In addition, we have $$[\alpha,\eta] = \lim\limits_{i \to \infty} [\alpha_i,\eta_i] = [\ti{\alpha},\ti{\eta}].$$ Thus, since $\ti{\eta} = u \eta$, we must have $\ti{\alpha} = \alpha \cdot \rec_{\Kz}(u)^{-1}$ in $\Gab_{\Kz} / N$. Putting all this together, we obtain
\bglnoz
  \lim\limits_{i \to \infty} [\psi(\alpha_i) \cdot \psi(\rec_{\Kz}(v_i)^{-1}),\xi_i] 
  &=& [\psi(\ti{\alpha}) \cdot \psi(\rec_{\Kz}(v)^{-1}),\xi]\\
  &=& [\psi(\alpha) \cdot \psi(\rec_{\Kz}(u)^{-1}) \cdot \psi(\rec_{\Kz}(v)^{-1}),\xi]\\
  &=& \Phi([\alpha,\eta]).
\eglnoz
This shows that $\Phi$ is continuous.

The same construction as for $\Phi$ gives rise to a continous map $\Psi: \ X_{\Lz} \to X_{\Kz}$ uniquely determined by $$\Psi([\beta,s_{\Lz}(\mfn)]) = [\psi^{-1}(\beta),s_{\Kz}(\varphi^{-1}(\mfn))]$$ for all $\beta \in \Gab_{\Lz}$ and $\mfn \in I_{\Lz}$. Obviously, $$\Phi \circ \Psi = \id_{X_{\Lz}} \mbox{ and }\Psi \circ \Phi = \id_{X_{\Kz}}.$$ Also, it is clear that $\Phi$ and $\Psi$ are monoid homomorphisms. Finally, for every $\mfm \in I_{\Kz}$, we have, because of \eqref{iv2}:
\bglnoz
  \Phi(\mfm) 
  &=& \Phi([\rec_{\Kz}(s_{\Kz}(\mfm))^{-1},s_{\Kz}(\mfm)] 
  = [\psi(\rec_{\Kz}(s_{\Kz}(\mfm)^{-1})),s_{\Lz}(\varphi(\mfm))]\\
  &=& [\rec_{\Lz}(s_{\Lz}(\varphi(\mfm)))^{-1},s_{\Lz}(\varphi(\mfm))] 
  = \varphi(\mfm).
\eglnoz
\eproof

\bremark
\label{reci-->monoid_N}
Note that in Proposition~\ref{reci-->monoid}, the monoid isomorphism $\varphi$ in (iv) and the conclusion is the same. In particular, if it was norm-preserving, then we will keep that property.
\eremark

\bprop[\cite{CLMSS}, 4.1]
\label{reci-->N}
Assume \textup{(iv)}.
Then
$  N(\mfp) = N(\varphi(\mfp))\ \text{for all} \ \mfp \in \cP_{\Kz}.$ \qed
\eprop

\section{Algebraic crossed products and orbit equivalence} \label{NC}

In this section, we study the relationship between the notions of equivalence used for the dynamical systems in Theorem \ref{THM} and isomorphisms of the corresponding algebraic crossed product algebras. We will prove our result in a more general context and then formulate the conclusion for our specific system. After that, we discuss the relation to quantum statistical mechanical systems such as the Bost-Connes system. 

Let $P \overset{\alpha}{\acts} X$ be an action (from the left) of a monoid $P$ on a locally compact Hausdorff space $X$. This defines a dynamical system $(X,P,\alpha)$. 
Two such dynamical systems $(X_i,P_i,\alpha_i)$ ($i=1,2$) are called \emph{orbit equivalent} if there exists a homeomorphism $\Phi: \ X_1 \to X_2$ with the property that for every $x, \, y \in X_1$, $x \in P_1 \cdot y$ $\LRarr$ $\Phi(x) \in P_2 \cdot \Phi(y)$. Here we write $P_1 \cdot y$ for $\menge{\alpha_p(y)}{p \in P_1}$.

Assume that $\alpha$ has the property that $\alpha^*_p(f) \defeq f \circ \alpha_p$ for $f \in C_0(X)$ defines a right action (called $\alpha^*$) of $P$ on $C_0(X)$ by endomorphisms of C*-algebras. (This is for example the case if all the $\alpha_p$ are proper.) From now on, all our dynamical systems are assumed to have this property. We form the algebraic crossed product $C_0(X) \rtimes^{{\rm alg}}_{\alpha^*} P$ of $(X,P,\alpha)$, as follows: As a complex vector space, $$C_0(X) \rtimes^{{\rm alg}}_{\alpha^*} P \cong \bigoplus_{p \in P} C_0(X).$$ A typical element of $C_0(X) \rtimes^{{\rm alg}}_{\alpha^*} P$ is of the form $\sum_p p f_p$. Here the sum is finite, and $p f_p$ stands for the vector in $\bigoplus_{p \in P} C_0(X)$ whose $p$-th coordinate is $f_p \in C_0(X)$ and whose remaining coordinates are zero. To define a structure of a complex algebra on $C_0(X) \rtimes^{{\rm alg}}_{\alpha^*} P$, we define $$(\sum_p p f_p) (\sum_q q g_q) \defeq \sum_{p,q} (pq) (\alpha^*_q(f_p) g_q).$$ In other words, we impose the commutation relation $f p = p \alpha^*_p(f)$ for $p \in P$ and $f \in C_0(X)$. Obviously, $C_0(X)$ embeds into $C_0(X) \rtimes^{{\rm alg}}_{\alpha^*} P$ as a subalgebra via $$C_0(X) \ni f \ma e_P f \in C_0(X) \rtimes^{{\rm alg}}_{\alpha^*} P,$$ where $e_P$ is the identity element of $P$.

The \emph{vanishing ideal} $I_x$ corresponding to $x \in X$ is the ideal $I_x = \menge{f \in C_0(X)}{f(x) = 0}$ in $C_0(X)$. 

\blemma
Let $(X,P,\alpha)$ be a system as above. For $x$ and $y$ in $X$, we have $$x = \alpha_p(y) \mbox{ if and only if }\alpha^*_p(I_x) \subseteq I_y.$$ 
\elemma
\bproof
\an{$\Rarr$}: If $x = \alpha_p(y)$, then every $f \in I_x$ satisfies $\alpha^*_p(f)(y) = f(\alpha_p(y)) = f(x) = 0$, and hence $\alpha^*_p(f) \in I_y$.

\an{$\Larr$}: To show that $x = \alpha_p(y)$, it suffices to show that every $f \in I_x$ satisfies $f(\alpha_p(y)) = 0$. But for arbitrary $f \in I_x$, $\alpha^*_p(I_x) \subseteq I_y$ implies that $\alpha^*_p(f) \in I_y$, so we get $f(\alpha_p(y)) = \alpha^*_p(f)(y) = 0$.
\eproof

Given a subspace $V$ of $C_0(X) \rtimes^{{\rm alg}}_{\alpha^*} P$ and $p \in P$, let $V_p \defeq \menge{f \in C_0(X)}{p f \in V}$.

\bprop
\label{xalpha_py}
Let $(X,P,\alpha)$ be a system as above, and let $A = C_0(X) \rtimes^{{\rm alg}}_{\alpha^*} P$. For $x$, $y$ in $X$ and $p \in P$, we have $$x = \alpha_p(y) \mbox{ if and only if }(\lspan(I_x A + A I_y))_p \neq C_0(X).$$
\eprop
\bproof
A straightforward computation shows
$$(\lspan(I_x A + A I_y))_p = \lspan(\alpha^*_p(I_x) C_0(X)) + I_y.$$
Since $I_y$ is a maximal (algebraic) ideal of $C_0(X)$, it follows that $$(\lspan(I_x A + A I_y))_p \neq C_0(X) \iff \lspan(\alpha^*_p(I_x) C_0(X)) \subseteq I_y \iff \alpha^*_p(I_x) \subseteq I_y.$$ By the previous lemma, the latter is equivalent to $x = \alpha_p(y)$.
\eproof

\bcor
\label{algcropro->oe}
Let $(X_i,P_i,\alpha_i)$ ($i=1,2$) be two systems as above. If there exists an isomorphism of complex algebras $C_0(X_1) \rtimes^{{\rm alg}}_{\alpha^*_1} P_1 \cong C_0(X_2) \rtimes^{{\rm alg}}_{\alpha^*_2} P_2$ which identifies the subalgebra $C_0(X_1)$ with $C_0(X_2)$, then $(X_1,P_1,\alpha_1)$ and $(X_2,P_2,\alpha_2)$ are orbit equivalent. \qed
\ecor

For a global field $\Kz$, let $C(X_{\Kz}) \rtimes^{{\rm alg}} I_{\Kz}$ be the algebraic crossed product corresponding to $I_{\Kz} \acts X_{\Kz}$.
\bcor
The statements in Theorem~\ref{THM} are equivalent to:
\begin{itemize}
\item[($*$)] There exists an isomorphism of $\Cz$-algebras $C(X_{\Kz}) \rtimes^{{\rm alg}} I_{\Kz} \overset{\cong}{\lori} C(X_{\Lz}) \rtimes^{{\rm alg}} I_{\Lz}$ which restricts to an isomorphism $C(X_{\Kz}) \overset{\cong}{\lori} C(X_{\Lz})$.
\end{itemize}
\ecor
\bproof
It follows from Corollary~\ref{algcropro->oe} that ($*$) implies Theorem~\ref{THM}~(i). Conversely, Theorem~\ref{THM}~(ii) clearly implies ($*$).
\eproof

\begin{remark} The (reduced) $C^*$-algebraic crossed product algebra $A_{\Kz}:=C(X_{\Kz}) \rtimes I_{\Kz}$ has a continuous one-parameter group of automorphism $$\sigma_{\Kz} \colon \Rz \rightarrow \Aut(A_{\Kz}) \colon t \mapsto \sigma_{\Kz,t} $$ given by extending to the entire algebra the map $$\sigma_{\Kz,t}(f)=f \mbox{ for } f \in C(X_{\Kz}) \mbox{ and }\sigma_{\Kz,t}(\mfp) = N(\mfp)^{it} \mbox{ for all }\mfp \in \cP_{\Kz}.$$ In this way, the pair $(A_{\Kz}, \sigma_{\Kz})$ becomes a \emph{quantum statistical mechanical system} in the sense of the first section of \cite{BC}. For $\Kz=\Qz$, this is the famous Bost-Connes system. For general number fields, it was introduced and studied in \cite{HP} and \cite{LLN}, and for function fields, in \cite{RN}. In \cite{CM1}, the first and last author proved that if two number fields $\Kz$ and $\Lz$ have isomorphic quantum statistical mechanical systems (in the obvious sense of isomorphism), then the number fields are arithmetically equivalent, i.e., they have the same Dedekind zeta function; essentially because this is their \emph{partition function}. By the theory of Ga{\ss}mann (\cite{Klingen}), this implies in particular that two such number fields that are Galois over $\Qz$ are isomorphic (without the need for any assumption on a subalgebra or of restricting to an algebraic crossed product, but making use of intertwining of the time evolution). It would be interesting to know in general whether an isomorphism of such systems $(A_{\Kz}, \sigma_{\Kz}) \cong (A_{\Lz}, \sigma_{\Lz})$ is the same as isomorphism of number fields $\Kz \cong \Lz$. 
\end{remark}

\end{document}